\documentclass[12pt]{article}
\newcommand{\copyleft}{
GNU FDL\thanks{
Copyright (C) 1976, 1985 Peter G. Doyle.
Permission is granted to copy, distribute and/or modify this document
under the terms of the GNU Free Documentation License, 
as published by the Free Software Foundation;
with no Invariant Sections, no Front-Cover Texts, and no Back-Cover Texts.
}}
\title{
A 27-vertex graph that is vertex-transitive  and edge-transitive
but not l-transitive
}
\author{Peter G. Doyle}
\date{Version dated 1985
\thanks{Derived from the Harvard senior thesis of Peter G. Doyle,
dated June 1976.}
\\ \copyleft
}

\usepackage{epsfig}
\newcommand{\fig}[2]{
\begin{figure}
\psfig{figure=figures/#1.ps,width=370pt}
\caption{#2}
\label{#1}
\end{figure}
}
\begin{document}
\maketitle

\begin{abstract}
I describe a 27-vertex graph that is vertex-transitive
and edge-transitive but not 1-transitive.
Thus while all vertices and edges of this graph are similar, there are no
edge-reversing automorphisms.
\end{abstract}
A graph (undirected, without loops or multiple edges) is said to be
{\em vertex-transitive} if its automorphism group acts transitively on the
set of vertices, {\em edge-transitive} if its automorphism group acts
transitively on the set of undirected edges, and {\em 1-transitive} if its
automorphism group acts transitively on the set of paths of length 1.
If a graph is edge-transitive but not 1-transitive then any edge can
be mapped to any other, but in only one of the two possible ways. In
my Harvard senior thesis
\cite{doyle:seniorThesis}, I described
a graph that is vertex-transitive and edge-transitive
but not 1-transitive. It has 27 vertices, and is regular of degree 4.
This beautiful graph was also discovered by Derek Holt
\cite{holt:doyleGraph}.
It seems likely
that this is the smallest graph that is vertex-transitive and edge-transitive
but not 1-transitive.

The question of the existence of graphs that are vertex-transitive
and edge-transitive but not 1-transitive was raised by Tutte
\cite{tutte:connectivity},
who showed
that any such graph must be regular of even degree. The first
examples were given by Bouwer
\cite{bouwer:graphs}.
Bouwer's smallest example has
54 vertices, and is regular of degree 4. While Bouwer's method of
construction differs from the method used here, Ronald Foster has
pointed out to me that the 27-vertex graph described here can be
obtained from Bouwer's 54-vertex graph by identifying pairs of
diametrically opposed vertices.

Recall that given a group $G$ and a set
$H \subseteq G - \{1\}$
such that $H=H^{-1}$, we
construct the group-graph $\Gamma_{G,H}$
by taking $G$ as the set of vertices,
and connecting every $g \in G$ to every element of the set $gH$. The idea,
which is inspired by work of Watkins
\cite{watkins:action},
will be to find a group $G$
and a set of generators $K \subseteq G-\{1\}$ such that:
\begin{enumerate}
\item
$K \cap K^{-1} = \emptyset$.
\item
For any $k_l ,k_2 \in K$, there is an automorphism $\phi$ of $G$
such that $\phi(k_1) = k_2$.
\item
If $\phi$ is an automorphism of $G$ such that
$\phi(K \cup K^{-1}) = K \cup K^{-1}$, then $\phi(K)=K$.
\end{enumerate}
The group-graph
$\Gamma_{G,K\cup K^{-1}}$
will then be vertex-transitive and edge-transitive,
and we may hope that conditions 1--3 will preclude its being 1-transitive.

For the group $G$ we take the non-abelian group of order 27 with
generators $a, b$ and relations
\[
a^9 = 1,\;\; b^3 = 1,\;\; b^{-1}ab = a^4
.
\]
(Cf. Hall
\cite{hall:groups}, p. 52.)
Setting $c= ba^{-1}$, we find that $G$ can be described
as the group with generators $a, c$ and relations
\[
a^9= 1,\;\; c^9= 1,
\]
\[
c^3= a^{-3},\;\; a^3= c^{-3},
\]
\[
c^{-1}ac=a^4,\;\; a^{-1}ca=c^4.
\]
These relations are not independent. Their redundancy allows us to
see at a glance that there is an automorphism $\phi$ of $G$ such that
$\phi(a)=c, \phi(c)=a$.
Setting $K= \{a,c\}$, we see that $K$ satisfies conditions 1 and
2, and it is easy to show that condition 3 also holds.

The graph
$\Gamma = \Gamma_{G,K\cup K^{-1}}$
is shown in Figure \ref{27}.
\fig{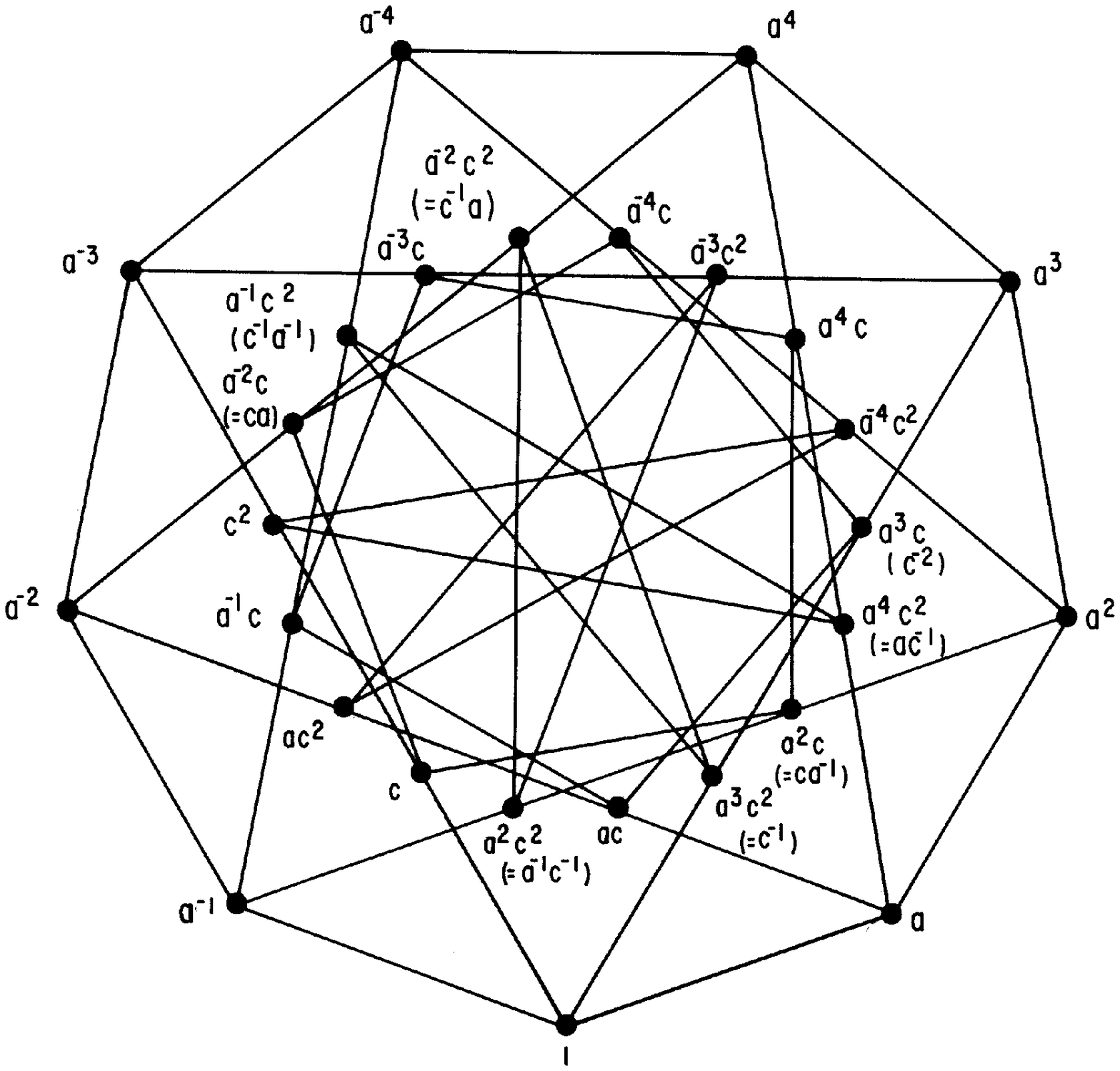}{The graph $\Gamma$.}
Although it is not obvious
from the drawing, we know that this graph
is vertex-transitive and edge-transitive. To
see that it is not 1-transitive, consider the subgraph $\Gamma^\prime$
obtained by
removing all vertices whose distance from the identity is $>2$.
(See Figure \ref{nbhd}.)
\fig{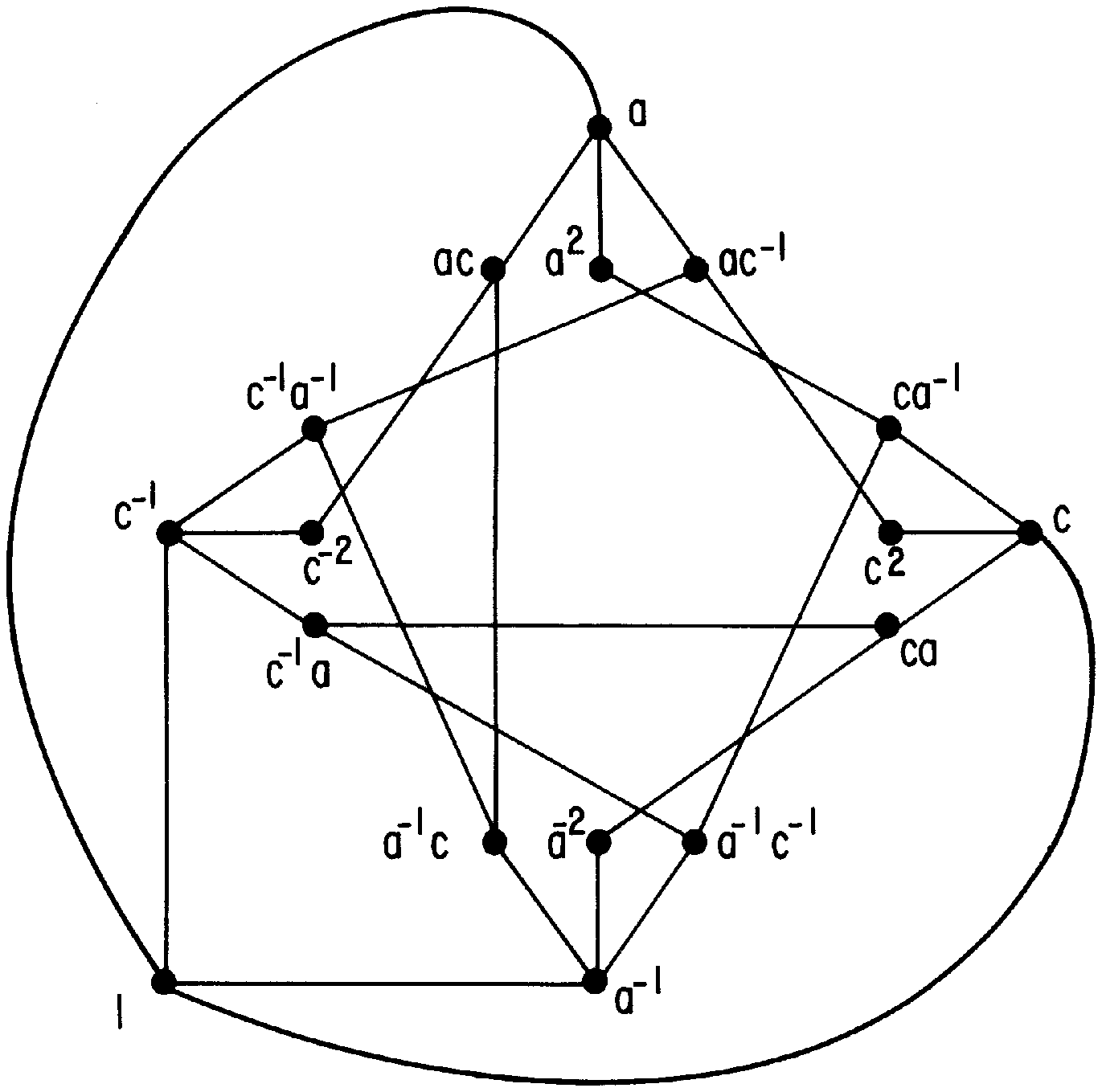}{The subgraph $\Gamma^\prime$.}
If there were a graph-automorphism $\phi$ of $\Gamma$ such that 
$\phi(1)=1, \phi(a)=a^{-1}$,
the restriction $\phi^\prime$ of $\phi$ to $\Gamma^\prime$
would be an automorphism of
$\Gamma^\prime$ such that
$\phi^\prime(1)=1,\phi^\prime(a)=a^{-1}$,
but it is easy to verify that no such
automorphism exists. Hence $\Gamma$ is not 1-transitive.

{\bf Acknowledgements.}
I would like to thank W. T. Tutte and Ronald Foster for helpful
correspondence.

\bibliography{bouwer}

\begin{thebibliography}{1}

\bibitem{bouwer:graphs}
I.~Z. Bouwer.
\newblock Vertex and edge transitive, but not 1-transitive graphs.
\newblock {\em Canadian Math. Bull.}, 13:231--237, 1970.

\bibitem{doyle:seniorThesis}
P.~G. Doyle.
\newblock On transitive graphs.
\newblock Senior Thesis, Harvard College, April 1976.

\bibitem{hall:groups}
M.~Hall.
\newblock {\em The Theory of Groups}.
\newblock Macmillan, New York, 1959.

\bibitem{holt:doyleGraph}
D.~F. Holt.
\newblock A graph which is edge transitive but not arc transitive.
\newblock {\em J. Graph Theory}, 5:201--204, 1981.

\bibitem{tutte:connectivity}
W.~T. Tutte.
\newblock {\em Connectivity in Graphs}.
\newblock University of Toronto Press, Toronto, 1966.

\bibitem{watkins:action}
M.~E. Watkins.
\newblock On the action of non-abelian groups on graphs.
\newblock {\em J. Combin. Theory}, 11:95--104, 1971.

\end{thebibliography}
\bibliographystyle{plain}

\end{document}